\documentclass[a5paper, 10pt]{article}

\usepackage{cite, amsmath, amssymb, booktabs, lastpage, graphicx}
\usepackage[hidelinks]{hyperref}

\usepackage{geometry}
\usepackage{setspace}

\usepackage[dvipsnames]{xcolor}
\usepackage{amsthm}

\theoremstyle{plain}

\theoremstyle{remark}

\theoremstyle{definition}

\usepackage{graphicx}
\makeatletter
\renewcommand{\maketitle}{
	\begin{center}
		\baselineskip=0.30in
		{\Large\bfseries \@title} \par
		\vspace{5mm}
		\baselineskip=0.2in
		{\large\bfseries \@author}\par
		\vspace{1mm}
		{\it \@address} \par
		{\small\tt \@email} \par
		\vspace{3mm}
		{\small (Received \@date)} \par
	\end{center}
	\vspace{3mm}
}
\newcommand{\address}[1]{\def\@address{#1}}
\newcommand{\email}[1]{\def\@email{#1}}

\makeatother

\newcommand{\acknowledgment}[1]{\vspace{5mm}\singlespacing
	{\noindent\textbf{\textit{Acknowledgment\/}:} #1}
}

\usepackage[margin=1cm,%
			font=footnotesize,%
			format=hang,%
			labelsep=period,%
			labelfont=bf]{caption}

\newcommand{\tbox}[1]{\mbox{\tiny #1}}
\newcommand{\bra}{\left\langle}
\newcommand{\ket}{\right\rangle}

\title{On the distribution of topological and spectral indices on random graphs}
\author{C. T. Mart\'inez-Mart\'inez$^{a,b,}$\footnote{Corresponding author.}, R. Aguilar-S\'anchez$^c$, \hspace{1cm}   J. A. M\'endez-Berm\'udez$^d$}

\address{
	$^a$Universidad Aut\'onoma de Guerrero, Centro Acapulco CP 39610, Acapulco de Ju\'arez, Guerrero, Mexico\\
    $^b$Facultad de Ciencias, Universidad Aut\'onoma Benito Ju\'arez de Oaxaca, Oaxaca de Ju\'arez 68120, Mexico\\
	$^c$Facultad de Ciencias Qu\'imicas, Benem\'erita Universidad Aut\'onoma de Puebla, Puebla 72570, Mexico\\
	$^d$Instituto de F\'{\i}sica, Benem\'erita Universidad Aut\'onoma de Puebla, Puebla 72570, Mexico
    }

\email{cl4ud7@gmail.com, ras747698@gmail.com, jmendezb@ifuap.buap.mx}

\date{\today}

\begin{document}

\maketitle

\begin{abstract}

We perform a detailed statistical study of the distribution of topological and spectral indices on random graphs 
$G=(V,E)$ in a wide range of connectivity regimes.
First, we consider degree-based topological indices (TIs), and focus on two classes of them:
$X_\Sigma(G) = \sum_{uv \in E} f(d_u,d_v)$ and 
$X_\Pi(G) = \prod_{uv \in E} g(d_u,d_v)$,
where $uv$ denotes the edge of $G$ connecting the vertices $u$ and $v$, $d_u$ is the degree of the vertex 
$u$, and $f(x,y)$ and $g(x,y)$ are functions of the vertex degrees.
Specifically, we apply $X_\Sigma(G)$ and $X_\Pi(G)$ on Erd\"os-R\'enyi graphs and random geometric graphs 
along the full transition from almost isolated vertices to mostly connected graphs.
While we verify that $P(X_\Sigma(G))$ converges to a standard normal distribution, we show that $P( X_\Pi(G))$ 
converges to a log-normal distribution.
In addition we also analyze Revan-degree-based indices and spectral indices (those defined from
the eigenvalues and eigenvectors of the graph adjacency matrix).
Indeed, for Revan-degree indices, we obtain results equivalent to those for standard degree-based TIs.
Instead, for spectral indices, we report two distinct patterns: 
the distribution of indices defined only from eigenvalues approaches a normal distribution, 
while the distribution of those indices involving both eigenvalues and eigenvectors approaches 
a log-normal distribution.

\end{abstract}

\onehalfspacing

\section{Introduction}

In chemical graph theory, graph invariants are widely used to characterize the structural properties of 
graphs~\cite{T18, CT08, R97}. These invariants can be classified mainly into two types: Topological 
indices (TIs) and multiplicative topological indices (MTIs). TIs are typically defined as sums over vertex 
or edge functions, such as
\begin{equation}
\begin{aligned}
& X_{\Sigma,F_V}(G) = \sum_{u \in V(G)} F_V(d_u) 
\\
& \qquad \qquad \qquad \mbox{or} 
\\
& X_{\Sigma,F_E}(G) = \sum_{uv \in E(G)} F_E(d_u,d_v)
\end{aligned}
\label{TI}
\end{equation}
while MTIs are defined as products over vertex or edge functions, such as
\begin{equation}
\begin{aligned}
& X_{\Pi,F_V}(G) = \prod_{u \in V(G)} F_V(d_u) 
\\
& \qquad \qquad \qquad \mbox{or} 
\\
& X_{\Pi,F_E}(G) = \prod_{uv \in E(G)} F_E(d_u,d_v) \, .
\end{aligned}
\label{MTI}
\end{equation}
Here $uv$ denotes the edge of the graph $G=(V(G),E(G))$ connecting the vertices $u$ and $v$,
$d_u$ is the degree of the vertex $u$, and $F_X(x)$ and $F_X(x,y)$ are appropriate chosen functions,
see e.g.~\cite{G13}.
While both $X_\Sigma(G)$ and $X_\Pi(G)$ are referred to as topological indices in the literature, to make
a distinction between them, here we name $X_\Sigma(G)$ as topological indices (TIs) and $X_\Pi(G)$ 
as {\it multiplicative} topological indices (MTIs).
Some prominent examples of TIs are the Randi\'c index~\cite{R75}, the Zagreb indices~\cite{G72}, and 
the Sombor index~\cite{G21}, while among the MTIs we can mention the Narumi-Katayama index~\cite{NK84} 
and the multiplicative versions of the Zagreb indices~\cite{TC10}. All these indices (to be defined later), 
among others, will be analyzed below.

More recently, a new class of TIs based on the Revan vertex degree has been 
introduced, see e.g.~\cite{K17, K18, KG22, KMRS22, AHMS22}. 
The \textit{Revan vertex degree} of the vertex $u$ is defined as 
\begin{equation}
r_u = \Delta + \delta - d_u, 
\label{revan}
\end{equation}
where $\Delta$ and $\delta$ are the maximum and minimum degrees among the vertices of the graph $G$, respectively.
 Revan-degree indices, defined as
 \begin{equation}
\begin{aligned}
& RX_\Sigma(G) = \sum_{uv \in E(G)} F(r_u,r_v)  
\\
& \qquad \qquad \qquad \mbox{or} 
\\
& RX_\Pi(G) = \prod_{uv \in E(G)} F(r_u,r_v)
\end{aligned}
\label{RTI}
\end{equation}
are the Revan analogs of standard TIs and MTIs, respectively. That is, $RX_{\Sigma}(G)$ is the Revan version of 
$X_{\Sigma}(G)$ and $RX_{\Pi}(G)$ is the Revan version of $X_{\Pi}(G)$. In this work, we also explore the distributions of 
Revan-degree indices.

Additionally, spectral indices, defined in terms of the eigenvalues and eigenvectors of the graph adjacency matrix, 
have gained attention due to their ability to capture global graph properties avoiding the problem of degeneracy, 
present in standard TIs; see e.g.~\cite{RB19}.
Specifically, we compute the so-called Rodr\'iguez-Vel\'azquez indices~\cite{RB19,AMRS21} as well as the graph 
energy~\cite{LSG12,GR20} and the subgraph centrality~\cite{ER05} (to be defined later).

The use of topological and spectral indices on random graphs is relatively recent.
Moreover, since a given parameter set represents an infinite-size ensemble of random graphs, the computation of a 
graph invariant on a single graph may be irrelevant. In contrast, the computation of the average value of a graph 
invariant over a large ensemble of random graphs, all characterized by the same parameter set, may provide useful 
average information about the full ensemble. This statistical approach, well known in random matrix theory (RMT) 
studies, has been recently applied to random graphs and networks by means of topological and spectral 
indices, see e.g.~\cite{AHMS22,AMRS21,MMRS20,MMRS21,AMRS20,AMRS24,AHMS20}.
In fact, the average value of some topological indices have been shown to be equivalent to standard RMT 
measures~\cite{AMRS20,AMRS21}.

While most studies of topological and spectral indices on random graphs have been focused on the average values
of the indices, just a few have considered their probability distribution functions numerically~\cite{MMRS20,MMRS21} 
and analytically~\cite{Y24,Y25}. Specifically, on the one hand, in Refs.~\cite{MMRS20,MMRS21} the probability 
distribution functions of the Randi\'c index, the harmonic index, the sum-connectivity index, the modified Zagreb index, 
and the inverse degree index on Erd\"os-R\'enyi graphs were reported.
On the other hand, in Refs.~\cite{Y24,Y25} the the probability distribution functions of TIs on Erd\"os-R\'enyi graphs
and of the Randi\'c index on random geometric graphs were studied; see also the related Refs.~\cite{Y23,YZ23}.

Therefore, in this work to go a step forward in the direction addressed by Refs.~\cite{MMRS20,MMRS21,Y24,Y25}, 
we conduct a comprehensive statistical (numeirical) analysis of the probability distribution functions 
(from now on we will call them just distributions) of TIs, MTIs, Revan-degree 
indices, and spectral indices on two types of random graphs: Erd\"os-R\'enyi graphs (ERGs) and 
random geometric graphs (RGGs).

This paper is organized as follows.
In Sec.~\ref{models} we introduce the graph models and the parameter settings to be used in the numerical analysis.
In Sec.~\ref{results} we report the distributions of TIs, MTIs, Revan-degree 
indices, and spectral indices on ERGs. To avoid text saturation, the results corresponding to RGGs are reported in 
the Appendix.
Finally, our findings are summarized in Sec.~\ref{conclusions}.

\section{Graph models and parameter settings}
\label{models}

ERGs~\cite{SR51,ER59,ER60}, $G_{\tbox{ERG}}(n,p)$, are formed by $n$ vertices connected 
independently with probability $p \in [0,1]$.
While RGGs~\cite{DC02,P03}, $G_{\tbox{RGG}}(n,r)$, consist of $n$ vertices uniformly and 
independently distributed on the unit square, where an edge connects two vertices if their Euclidean 
distance is less or equal than the connection radius $r \in [0,\sqrt{2}]$.

Our study spans the full transition from almost isolated nodes ($p\to 0$ or $r\to 0$) to almost complete
graphs ($p\to 1$ or $r\to \sqrt{2}$), providing insights into the behavior of the indices across different 
connectivity regimes. In order to clearly set the connectivity regime on both random graph
models we will use the average number of non isolated vertices $\bra V(G) \ket$, which can also be
regarded as a TI, see e.g.~\cite{AMRS20}.

\begin{figure}
\centering
\includegraphics[width=0.7\textwidth]{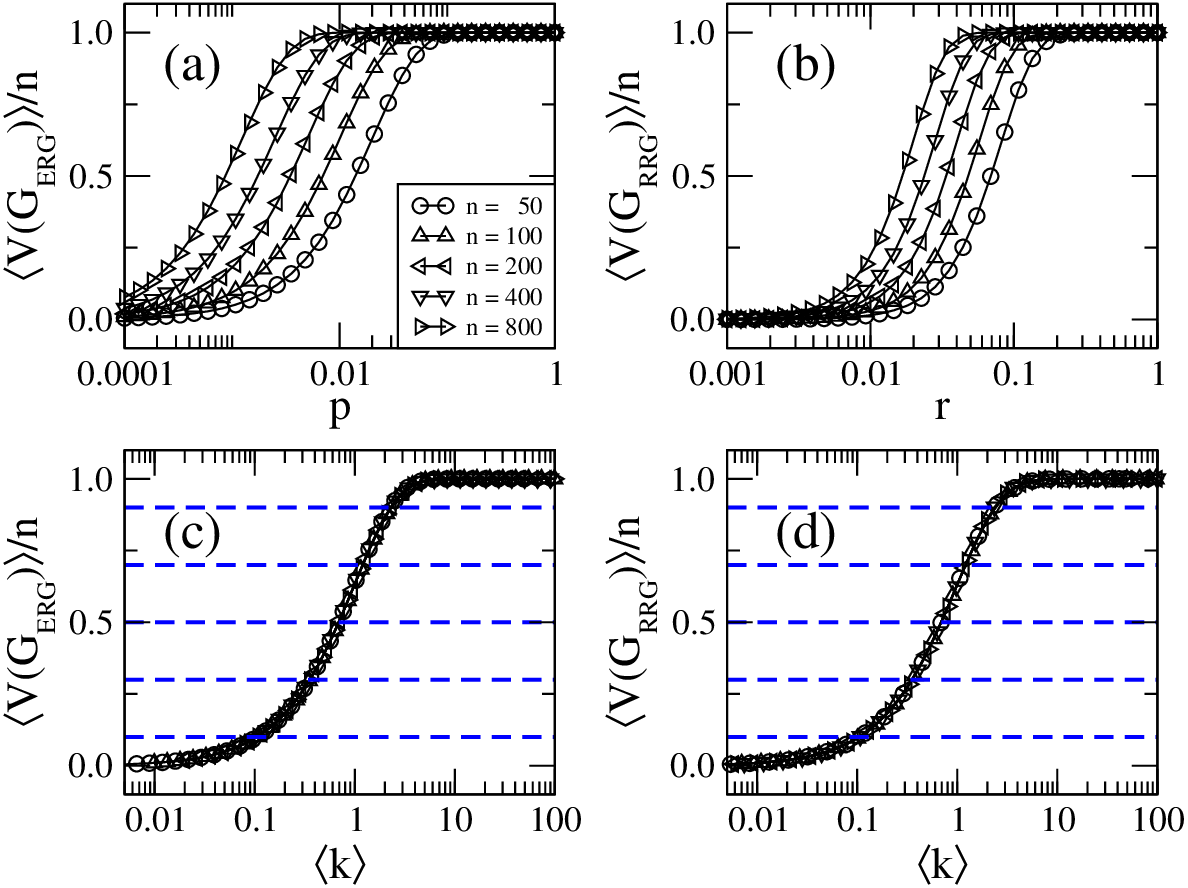}
\caption{Average number of non isolated vertices $\bra V(G) \ket$, normalized to the graph size $n$,
for Erd\H{o}s-R\'{e}nyi graphs as a function of (a) the probability $p$ and (c) the average degree 
$\bra k \ket = p(n-1)$.
$\bra V(G) \ket/n$ for random geometric graphs as a function of (b) the connection radius $r$
and (d) the average degree $\bra k \ket=g(r)(n-1)$, see Eq.~(\ref{gofr}).
The blue horizontal dashed lines in (c,d) indicate the values of $\bra V(G) \ket/n$ used to
construct the histograms in Figs.~\ref{Fig02}-\ref{Fig12}: $\bra V(G) \ket/n=0.1$, 0.3, 0.5, 0.7 and 0.9.
%The red full curves in (c,d) are Eq.~(\ref{Vofk}).
Each data value was computed by averaging over $10^{6}$ random graphs $G$.}
\label{Fig01}
\end{figure}

Since $\bra V(G) \ket = 0$ for graphs with only isolated nodes and $\bra V(G) \ket = n$ when all nodes are connected, $\bra V(G) \ket$ 
shows a smooth transition from $0$ to $n$ by increasing $p$ from 0 to 1 for ERGs or by increasing $r$ from 
0 to $\sqrt{2}$ for RGGs. This is indeed shown in Figs.~\ref{Fig01}(a) and~\ref{Fig01}(b) where we
present $\bra V(G) \ket$, normalized to the graph size $n$, for ERGs as a function of $p$ and for RGGs
as a function of $r$, respectively. There, different symbols correspond to different graph sizes. 

As well as for other TIs (see e.g.~\cite{AMRS24,MMRS20,AHMS20,MMRS21}) the average 
degree $\bra k \ket$ serves
as the scaling parameter of $\bra V(G) \ket/n$~\cite{AMRS20}; meaning that the curve $\bra V(G) \ket/n$ 
vs.~$\bra k \ket$ is a universal curve. This is verified in Figs.~\ref{Fig01}(c) and~\ref{Fig01}(d) where we
plot, respectively, $\bra V(G) \ket/n$ for ERGs and RGGs as a function of $\bra k \ket$: I.e., curves 
for different graph sizes fall one on top of the other. 
Note that the functional dependence of $\bra k \ket$ on the graph  
parameters is significantly different for both graph models; while for ERGs
\begin{equation}
\bra k \ket =  p(n-1) ,
\label{kERG}
\end{equation}
for RGGs it takes the form
\begin{equation}
\bra k \ket =  g(r)(n-1) 
\label{kRRG}
\end{equation}
with~\cite{EM15}
\begin{equation}
g(r) =
\left\{
\begin{array}{ll}
           r^2  \left[ \pi - \frac{8}{3}r +\frac{1}{2}r^2 \right] , & 0 \leq r \leq 1 \, ,
           \vspace{0.25cm} \\
             \frac{1}{3} - 2r^2 \left[ 1 - \arcsin(1/r) + \arccos(1/r) \right] \vspace{0.15cm} \\
             \qquad +\frac{4}{3}(2r^2+1) \sqrt{r^2-1} -\frac{1}{2}r^4 , & 1 \leq r \leq \sqrt{2} \, .
\end{array}
\right.
\label{gofr}
\end{equation}

Then, the curves $\bra V(G) \ket/n$ vs.~$\bra k \ket$ in Figs.~\ref{Fig01}(c) and~\ref{Fig01}(d)
allow us to set both graph models in a given connectivity regime regardless of the parameter
combinations. Specifically, in order to span the full transition from almost isolated nodes 
($\bra k \ket\to 0$) to almost complete graphs ($\bra k \ket\gg 1$), we choose five values of 
the ratio $\bra V(G) \ket/n$: 0.1, 0.3, 0.5, 0.7 and 0.9; as indicated by the blue horizontal 
dashed lines in in Figs.~\ref{Fig01}(c) and~\ref{Fig01}(d).

\begin{figure}
\centering
\includegraphics[width=0.95\textwidth]{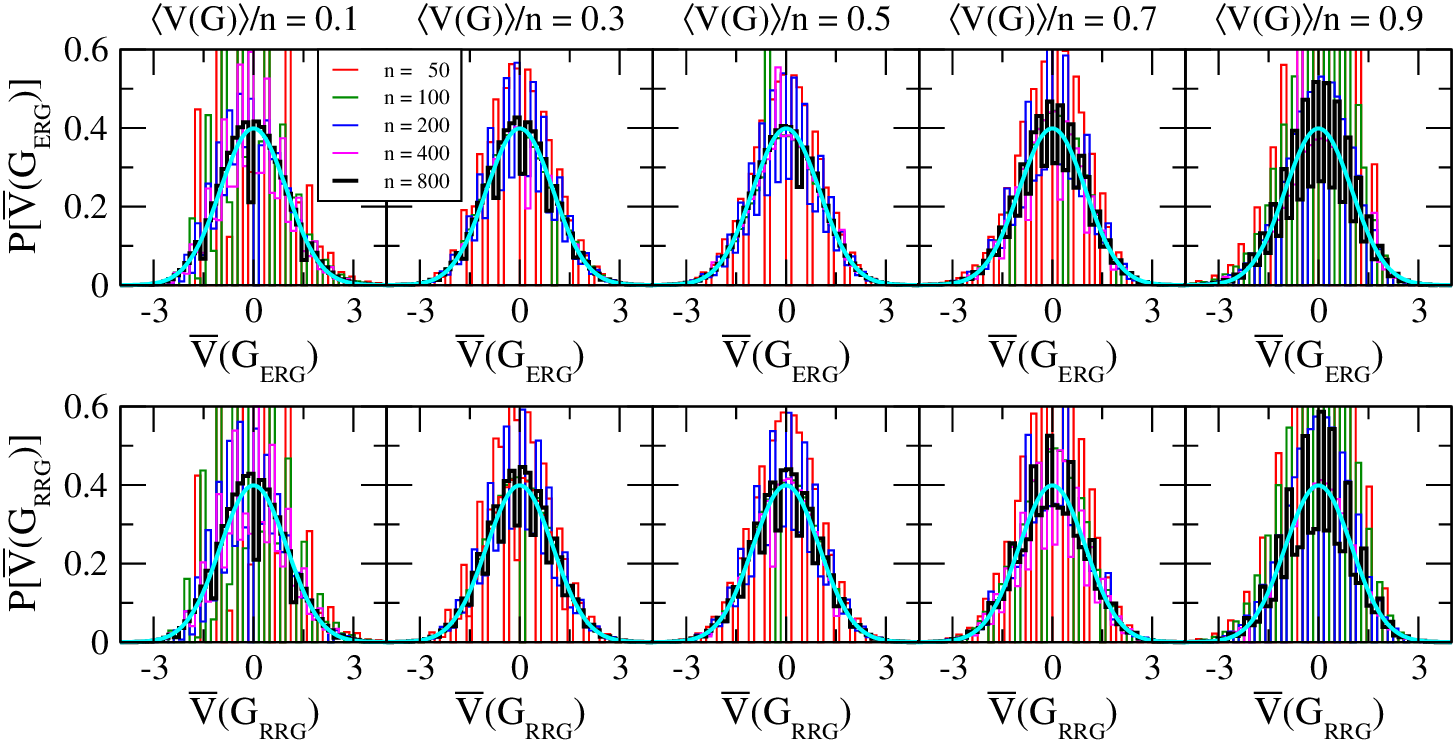}
\caption{Probability distribution functions of the normalized number non isolated vertices 
$\overline{V}(G)$ for
Erd\H{o}s-R\'{e}nyi graphs (upper panels) and for random geometric graphs (lower panels).
Each panel displays five histograms corresponding to graphs of different sizes $n$.
Each column corresponds to a fixed value of the ratio $\bra V(G) \ket/n$.
Each histogram is constructed with $10^{6}$ values of $V(G)$.
The cyan full line in all panels is a normal distribution with zero mean and unit variance.}
\label{Fig02}
\end{figure}

As a first example, in Fig.~\ref{Fig02} we present the probability distribution function of the 
normalized number of non isolated vertices $\overline{V}(G)=\bra V(G) \ket/n$ for both 
ERGs (upper panels) and RGGs (lower panels) for five values of $\bra V(G) \ket/n$
(0.1, 0.3, 0.5, 0.7 and 0.9); as indicated on top of the figure. Data were standardized to a zero mean and unit variance. 
Notice that each panel displays five histograms corresponding to graphs of different sizes $n$.
From this figure we can see, for both random graph models, that the distribution of this 
degree-based TI converges to a standard normal distribution (represented by the cyan full lines) 
as the graphs become larger regardless of the value of $\bra V(G) \ket/n$; see that black histograms 
in all panels approach the normal distribution.

\section{Distribution of topological and spectral indices on Erd\H{o}s-R\'{e}nyi graphs}
\label{results}

In this section, we analyze the distributions of TIs, MTIs, Revan-degree indices, 
and spectral indices on ERGs across the full range of connectivity.
Results corresponding to RGGs are reported in the Appendix.

\subsection{Distribution of degree-based topological indices on Erd\H{o}s-R\'{e}nyi graphs}
\label{Sec:TIs}

To explore the distribution of TIs, we selected the following well-known indices:
The first and second Zagreb indices~\cite{G72}
\begin{equation}
M_1(G) = \sum_{u\in V(G)} d_u^2 = \sum_{uv\in E(G)} d_u + d_v
\label{M1}
\end{equation}
and
\begin{equation}
M_2(G) = \sum_{uv\in E(G)} d_ud_v ,
\label{M2}
\end{equation}
respectively, the Sombor index~\cite{G21}
\begin{equation}
SO(G) = \sum_{uv\in E(G)} \sqrt{d_u^2 + d_v^2} ,
\label{SO}
\end{equation}
the Randi\'c connectivity index~\cite{R75}
\begin{equation}
R(G) = \sum_{uv\in E(G)} \frac{1}{\sqrt{d_ud_v}} ,
\label{R}
\end{equation}
and the harmonic index~\cite{F87}
\begin{equation}
H(G) = \sum_{uv\in E(G)} \frac{2}{d_u + d_v} .
\label{H}
\end{equation}

Based on the results from the previous section, regarding the number of non-isolated vertices $V(G)$, 
we explore the distribution of TIs in three representative connectivity regimes:

(i) {\it Sparse regime}. When most vertices are isolated, $\bra V(G) \ket/n=0.1$.

(ii) {\it Intermediate regime}. When the proportion of isolated and non-isolated vertices 
is approximately equal, $\bra V(G) \ket/n=0.5$.

(iii) {\it Dense regime}. When most vertices are connected, $\bra V(G) \ket/n=0.9$.

\begin{figure}
\centering
\includegraphics[width=0.75\textwidth]{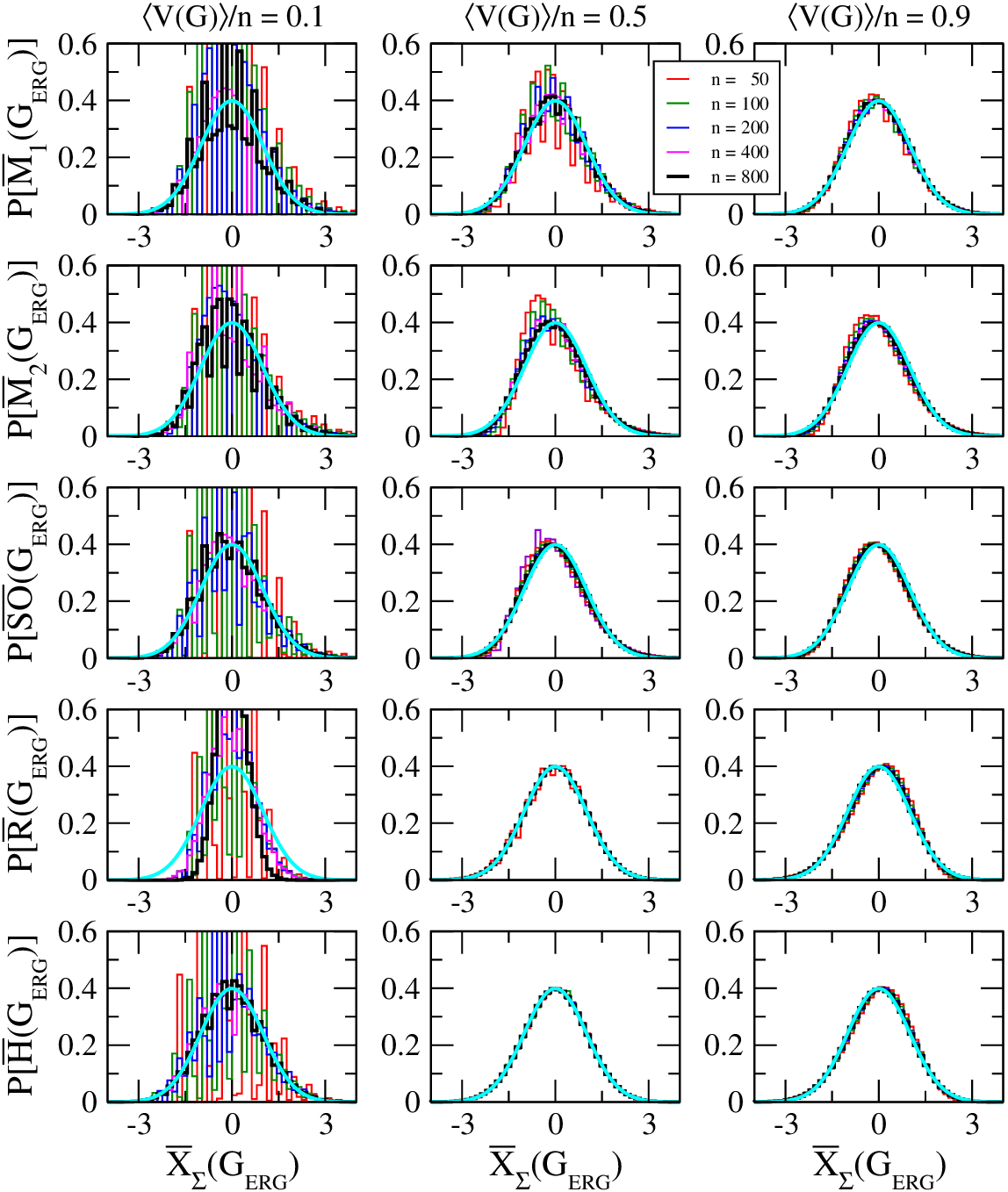}
\caption{Probability distribution functions of standardized degree-based topological indices on 
Erd\H{o}s-R\'{e}nyi graphs: 
First Zagreb index $\overline{M}_1(G)$, second Zagreb index $\overline{M}_2(G)$, 
Sombor index $\overline{SO}(G)$, Randi\'c index $\overline{R}(G)$, and harmonic index 
$\overline{H}(G)$.
Each panel displays five histograms corresponding to graphs of different sizes $n \in [50, 800]$.
Each column corresponds to a fixed value of the ratio $\bra V(G) \ket/n$.
Each histogram is constructed with $10^6$ values of $X_\Sigma(G)$. 
The cyan full line in all panels is a normal distribution with zero mean and unit variance.}
\label{Fig03}
\end{figure}

In Fig.~\ref{Fig03} we present the probability distribution functions of the TIs of Eqs.~(\ref{M1})-(\ref{H})
on ERGs. 
In this and all the following figures, each panel displays five histograms corresponding to graphs of different 
sizes $n \in [50, 800]$. Each histogram is constructed from an ensemble of $10^{6}$ random graphs.
Also, each column corresponds to a fixed value of the ratio $\bra V(G) \ket/n$; so, graph connectivity increases 
from left to right.
Moreover, to ease the comparison across regimes and graph sizes, the data is standardized to have zero 
mean and unit variance:
\begin{equation}
 \overline{X}_{\Sigma}(G)= \frac{ X_{\Sigma}(G)-\mu_{X_{\Sigma}}}{\sigma_{X_\Sigma}},
\end{equation}
where $\mu_{X_{\Sigma}}$ and $\sigma_{X_\Sigma}$ denote, respectively,  the mean and standard deviation 
of the TI $X_{\Sigma}(G)$.

From this figure, we can clearly observe that the distribution of all the TIs analyzed here tends to a normal 
distribution (indicated with the cyan line in all panels). 
It is interesting to note that the normal distribution is approached even in the sparse regime for all TIs 
(except for the Randi\'c index) when the graph size is large enough; see the black histograms in left panels
corresponding to $\bra V(G) \ket/n=0.1$.
Evidently, for $\bra V(G) \ket/n=0.5$ and 0.9 all histograms, even those corresponding to $n=50$, for all TIs
fall on top of the normal distribution.
Note that with Fig.~\ref{Fig03} we numerically validate the analytical results of Ref.~\cite{Y23} were the 
distribution of TIs was predicted to converge to a normal distribution.

We now proceed to compute the distributions of MTIs. To this end we consider the following well-known MTIs:
The Narumi-Katayama index~\cite{NK84}
\begin{equation}
NK(G) = \prod_{u\in V(G)} d_u ,
\label{NK}
\end{equation}
multiplicative versions of the Zagreb indices~\cite{TC10}:
\begin{equation}
\Pi_1(G) = \prod_{u\in V(G)} d_u^2 ,
\label{P1}
\end{equation}
\begin{equation}
\Pi_2(G) = \prod_{uv\in E(G)} d_ud_v ,
\label{P2}
\end{equation}
and
\begin{equation}
\Pi_1^*(G) = \prod_{uv\in E(G)} d_u + d_v ,
\label{P1*}
\end{equation}
the multiplicative Randi\'c connectivity index~\cite{AMRS24} 
\begin{equation}
R_\Pi(G) = \prod_{uv\in E(G)} \frac{1}{\sqrt{d_ud_v}} , 
\label{mR}
\end{equation}
and the multiplicative harmonic index~\cite{AMRS24}
\begin{equation}
H_\Pi(G) = \prod_{uv\in E(G)} \frac{2}{d_u + d_v} .
\label{mH}
\end{equation}

\begin{figure}
\centering
\includegraphics[width=0.75\textwidth]{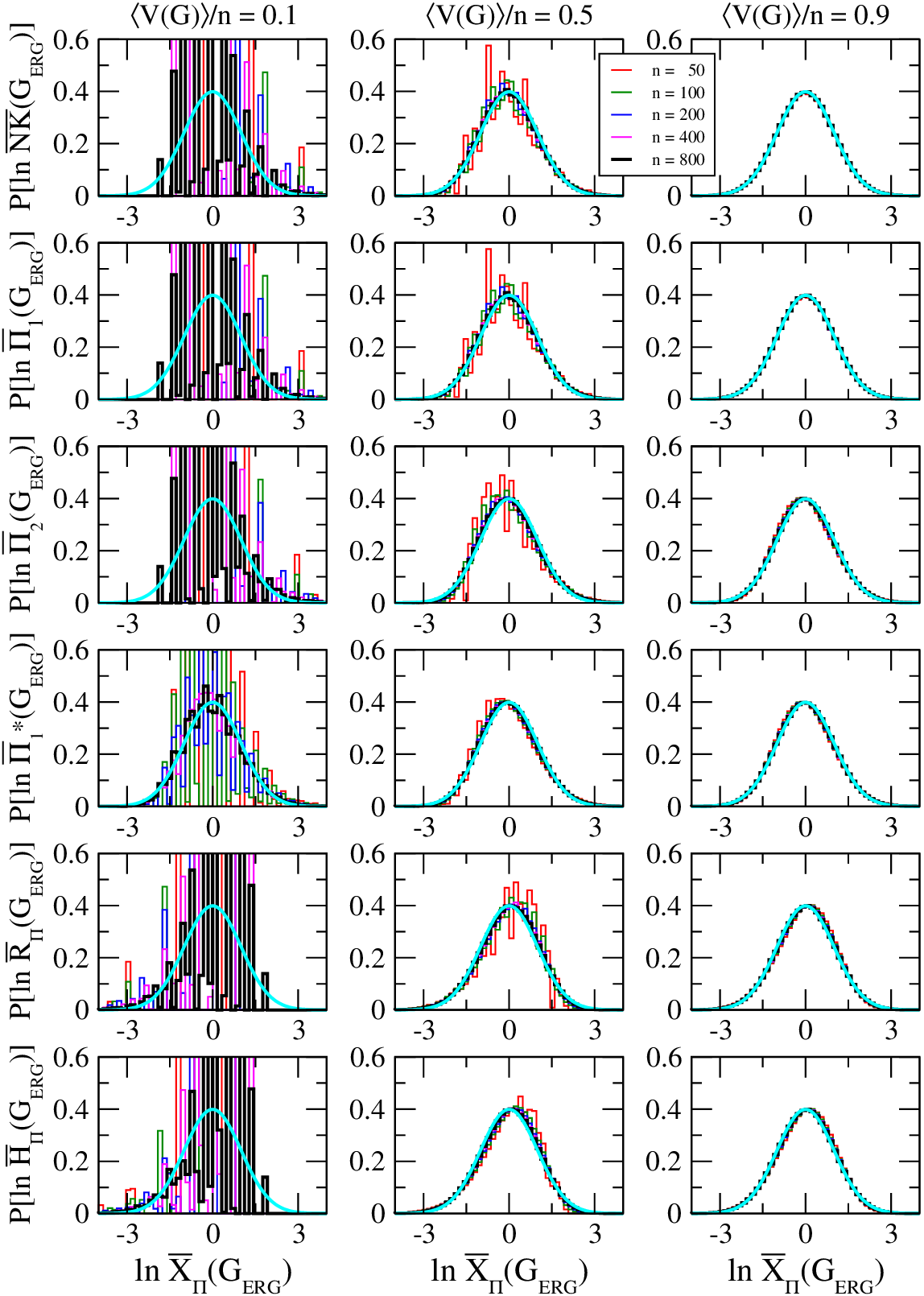}
\caption{Probability distribution functions of the logarithm of standardized multiplicative 
topological indices on Erd\H{o}s-R\'{e}nyi graphs: 
Narumi-Katayama index $\overline{NK}(G)$, multiplicative Zagreb indices $\overline{\Pi}_1(G)$, 
$\overline{\Pi}_2(G)$ and $\overline{\Pi}_1^*(G)$, multiplicative Randi\'c index $\overline{R}_\Pi(G)$, 
and multiplicative harmonic index $\overline{H}_\Pi(G)$.
Each panel displays five histograms corresponding to graphs of different sizes $n \in [50, 800]$.
Each column corresponds to a fixed value of the ratio $\bra V(G) \ket/n$.
Each histogram is constructed with $10^{6}$ values of $X_\Pi(G)$.
The cyan full line in all panels is a normal distribution with zero mean and unit variance.}
\label{Fig04}
\end{figure}

Since the values of MTIs grow exponentially for increasing average degree, see e.g.~\cite{AMRS24},
we compute their logarithms instead. 
Then, in Fig.~\ref{Fig04} we present the probability distribution function of the logarithm of the standardized MTIs
of Eqs.~\eqref{NK}-\eqref{mH} on ERGs.
From this figure, we can clearly observe that the distribution of the logarithm all the MTIs analyzed here tends 
to a normal distribution (indicated with the cyan line in all panels). 
Remarkably, the normal distribution is approached even in the sparse regime for all MTIs  
(not shown here) when the graph size is large enough; in fact, $n=800$ (the larger graph size we used in this work)
is not enough to observe the normal distribution in the sparse regime (see the black histograms in left panels
corresponding to $\bra V(G) \ket/n=0.1$).
For $\bra V(G) \ket/n=0.5$ and 0.9 all histograms, even those corresponding to $n=50$, for all MTIs
fall on top of the normal distribution.

Therefore, since the distributions of the logarithm of the MTIs follow a normal distribution, we can conclude 
that the distributions of the MTIs follow a log-normal distribution.

\subsection{Distribution of Revan-degree indices on Erd\H{o}s-R\'{e}nyi graphs}

The Revan-degree indices we choose for our study are: 
The first and second Revan Zagreb indices~\cite{K17}
\begin{equation}
R_1(G) = \sum_{u\in V(G)} r_u^2 = \sum_{uv\in E(G)} r_u + r_v
\label{rM1}
\end{equation}
and
\begin{equation}
R_2(G) = \sum_{uv\in E(G)} r_ur_v ,
\label{rM2}
\end{equation}
respectively, 
the Revan Sombor index~\cite{KG22}
\begin{equation}
RSO(G) = \sum_{uv\in E(G)} \sqrt{r_u^2 + r_v^2} ,
\label{rSO}
\end{equation}
the Revan Randi\'c index
\begin{equation}
RR(G) = \sum_{uv\in E(G)} \frac{1}{\sqrt{r_ur_v}} ,
\label{rR}
\end{equation}
and the Revan harmonic index
\begin{equation}
RH(G) = \sum_{uv\in E(G)} \frac{2}{r_u + r_v} .
\label{rH}
\end{equation}
We note that, as far as we know, $RR(G)$ and $RH(G)$ are being introduced here.

We also explore the behavior of the distribution of the multiplicative Revan-degree indices, $RX_\Pi(G)$:
The multiplicative Revan Narumi-Katayama index
\begin{equation}
RNK(G) = \prod_{u\in V(G)} r_u ,
\label{rNK}
\end{equation}
the multiplicative Revan Zagreb indices
\begin{equation}
R_{1\Pi}(G) = \prod_{u\in V(G)} r_u^2 ,
\label{rP1a}
\end{equation}
\begin{equation}
R_{1\Pi^{*}}(G) = \prod_{uv\in E(G)} r_u+r_v ,
\label{rP1}
\end{equation}
and
\begin{equation}
R_{2\Pi}(G) = \prod_{uv\in E(G)} r_ur_v ,
\label{rP2}
\end{equation}
the multiplicative Revan Randi\'c connectivity index
\begin{equation}
RR_\Pi(G) = \prod_{uv\in E(G)} \frac{1}{\sqrt{r_ur_v}} , 
\label{rmR}
\end{equation}
and the multiplicative Revan harmonic index
\begin{equation}
RH_\Pi(G) = \prod_{uv\in E(G)} \frac{2}{r_u + r_v} .
\label{rmH}
\end{equation}
It is fair to mention that $R_{1\Pi}(G)$ and $R_{2\Pi}(G)$ were introduced in~\cite{AHMS22}; while,
as far as we know, $RNK(G)$, $R_{1\Pi^{*}}(G)$, $RR_\Pi(G)$ and $RH_\Pi(G)$ are being introduced here.

\begin{figure}
\centering
\includegraphics[width=0.75\textwidth]{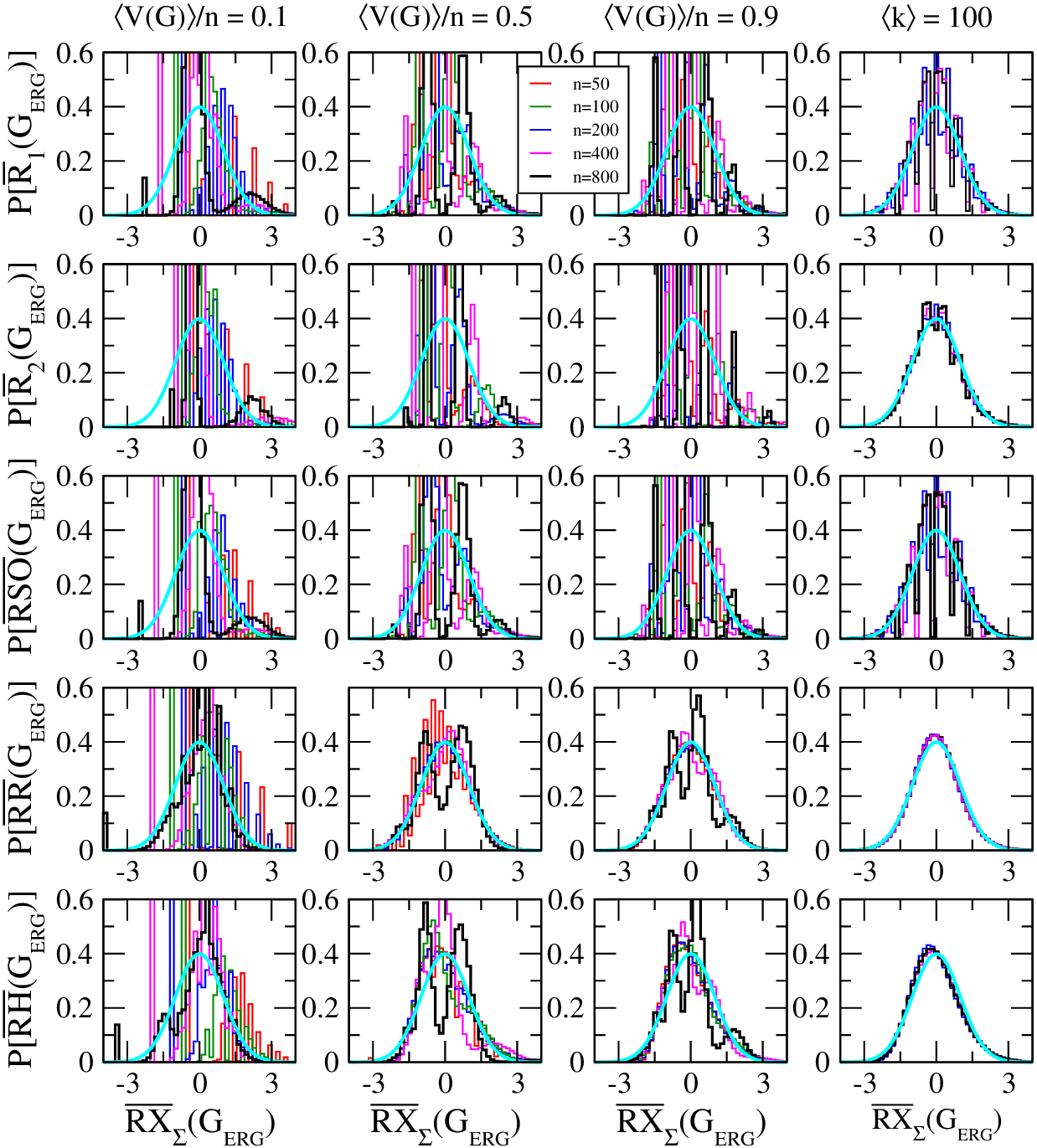}
\caption{Probability distribution functions of standardized Revan-degree indices on 
Erd\H{o}s-R\'{e}nyi graphs: 
First Revan Zagreb index $\overline{R}_1(G)$, second Revan Zagreb index $\overline{R}_2(G)$, 
Revan Sombor index $\overline{RSO}(G)$, Revan Randi\'c index $\overline{RR}(G)$, and 
Revan harmonic index $\overline{RH}(G)$.
Each panel displays five histograms corresponding to graphs of different sizes $n \in [50, 800]$.
Each column corresponds to a fixed value of the ratio $\bra V(G) \ket/n$, except for the right 
column where $\bra k \ket=100$ is set.
Each histogram is constructed with $10^{6}$ values of $RX_\Sigma(G)$.
The cyan full line in all panels is a normal distribution with zero mean and unit variance.}
\label{Fig07}
\end{figure}

\begin{figure}
\centering
\includegraphics[width=0.75\textwidth]{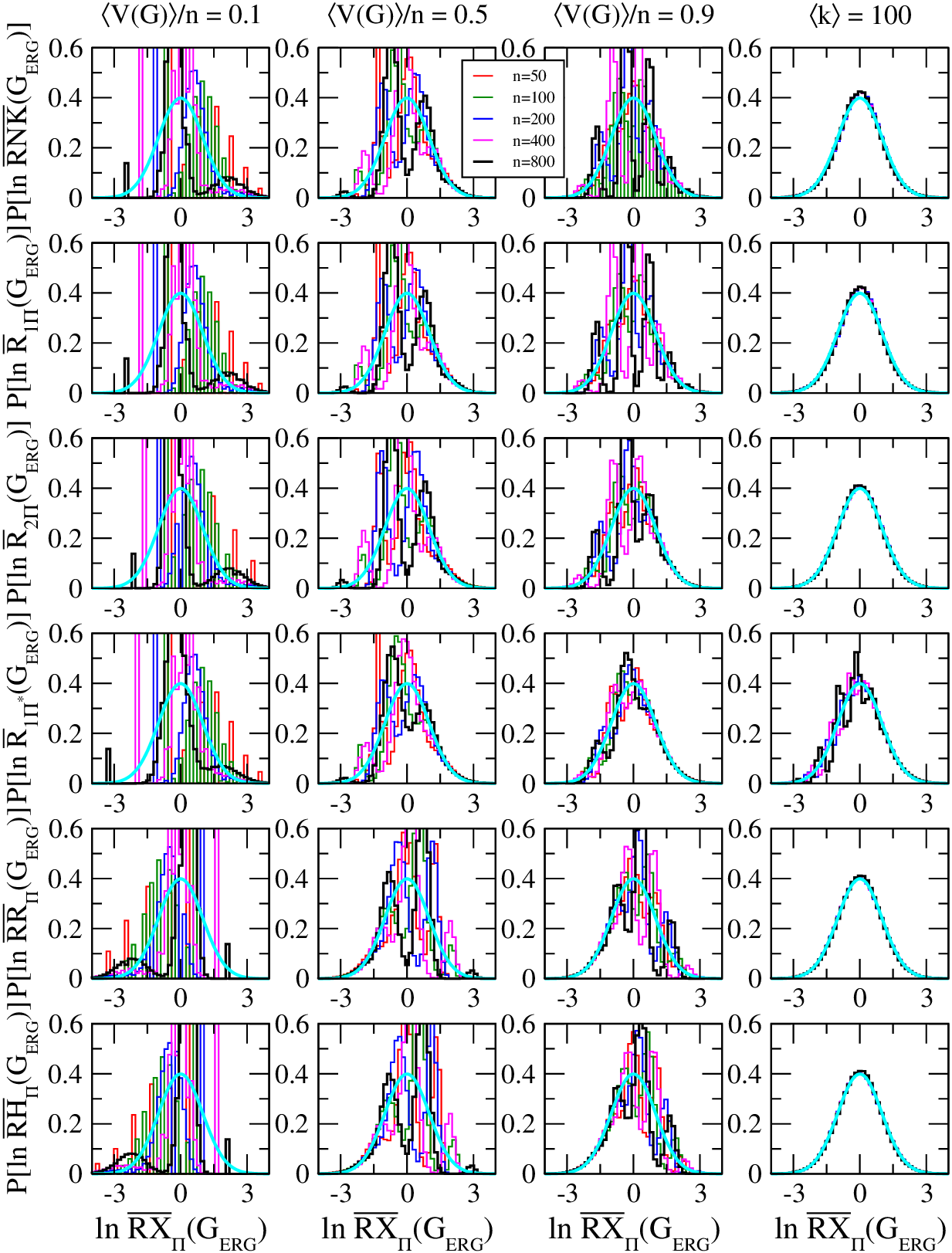}
\caption{Probability distribution functions of the logarithm of standardized multiplicative Revan-degree
indices on Erd\H{o}s-R\'{e}nyi graphs: 
Revan Narumi-Katayama index $\overline{RNK}(G)$, multiplicative Revan Zagreb indices $\overline{R}_{1\Pi}(G)$, 
$\overline{R}_{2\Pi}(G)$ and $\overline{R}_{1\Pi^*}(G)$,
multiplicative Revan Randi\'c index $\overline{RR}_\Pi(G)$, and multiplicative Revan harmonic index $\overline{RH}_\Pi(G)$.
Each panel displays five histograms corresponding to graphs of different sizes $n \in [50, 800]$.
Each column corresponds to a fixed value of the ratio $\bra V(G) \ket/n$, except for the right 
column where $\bra k \ket=100$ is set.
Each histogram is constructed with $10^{6}$ values of $RX_\Pi(G)$.
The cyan full line in all panels is a normal distribution with zero mean and unit variance.}
\label{Fig08}
\end{figure}

In Figs.~\ref{Fig07} and~\ref{Fig08} we present, respectively, the probability distribution functions of the
standardized Revan-degree indices of Eqs.~(\ref{rM1})-(\ref{rH}) and the probability distribution functions of 
the logarithm of the standardized multiplicative Revan-degree indices of Eqs.~\eqref{rNK}-\eqref{rmH},
both on ERGs. In contrast with TIs and MTIs, see Figs.~\ref{Fig03} and~\ref{Fig04}, 
we do not observe a clear transition of the distributions of 
Revan-degree indices nor of the distributions of the logarithm of multiplicative Revan-degree indices 
to a normal distribution; not even in the dense regime (see the panels in the third columns of 
Figs.~\ref{Fig07} and~\ref{Fig08} corresponding to $\bra V(G) \ket/n=0.9$).
However, in Ref.~\cite{AHMS22} it was shown that the statistical properties of both Revan-degree indices
and multiplicative Revan-degree indices are equivalent to those of their standard-degree counterparts
in the dense limit, specifically for $\bra k \ket > 10$.
Thus, we indeed expect to recover normal distributions of both Revan-degree indices the logarithm of 
multiplicative Revan-degree indices deep enough in the dense limit.


Therefore, we include an additional column in both Figs.~\ref{Fig07} and~\ref{Fig08} where we set 
$\bra k \ket$ to 100; i.e.~the graphs are now deeper in the dense regime.
So, we can now conclude that, deep in the dense regime, the distributions of 
Revan-degree indices follow normal distributions while the distributions of multiplicative Revan-degree indices
follow log-normal distributions.

\subsection{Distribution of spectral indices on Erd\H{o}s-R\'{e}nyi graphs}

Finally, we extend our analysis to spectral indices. To this end, we first define a weighted adjacency matrix as follows:
\begin{equation}
A_{ij}=\left\{
\begin{array}{cl}
\sqrt{2} \epsilon_{ii} \ & \mbox{for $i=j$}, \\
\epsilon_{ij} & \mbox{if there is an edge between vertices $i$ and $j$},\\
0 \ & \mbox{otherwise}.
\end{array}
\right.
\label{Aij}
\end{equation}

Here, $\epsilon_{ij}$ are statistically independent random variables drawn from a normal distribution with zero mean and unit variance. Once the adjacency matrix is weighted, we diagonalize it and compute the corresponding spectral indices. The spectral indices we consider include 
Rodr\'iguez-Vel\'azquez indices, the graph energy, and centrality-based indices, which are defined as 
follows.

Given an orthonormal basis of eigenvectors $\{\Psi_{i}\}_{i=1}^{n}$  and the corresponding eigenvalues 
$\{\lambda_{i}\}_{i=1}^{n}$ of the adjacency matrix $A$ of a graph $G$ of size $n$, the first and second 
Rodr\'iguez-Vel\'azquez (RV) indices are defined as~\cite{RB19}
\begin{equation}
RV_{a}(G)=\left(\frac{1}{n}\sum_{i=1}^{n}S_{i}^2\right)^{1/2}
\end{equation}
and
\begin{equation}
RV_{b}(G)=\sum_{i=1}^{n}x_{i}S_{i},
\end{equation}
respectively. Here,
\begin{equation} 
S_{i}=\sum_{j=1}^{n} \left(\Psi^{i}_{j}\right)^2 \exp(\lambda_{j})
\end{equation}
represents the subgraph centrality while
\begin{equation} 
x_{i}=\frac{1}{\lambda_{1}}\sum_{j=1}^{n}A_{ij}|\Psi_{j}^{1}| 
\end{equation}
denotes the eigenvector centrality of vertex $i$, where $\lambda_1$ is the largest eigenvalue of $A$ 
and $\Psi_{j}^{1}$ is the $jth$ component of the eigenvector corresponding  to $\lambda_1$. 

Additionally, the graph energy $E(G)$~\cite{LSG12,GR20} and the exponential subgraph centrality $EE(G)$~\cite{ER05} 
are defined as
\begin{equation}
E(G)=\sum_{i=1}^{n} |\lambda_{i}|,
\end{equation}
and
\begin{equation}
EE(G)=\frac{1}{n}\sum_{i=1}^{n} S_{i}=\frac{1}{n}\sum_{i=1}^{n} \exp{(\lambda_{i})},
\end{equation}
respectively.

\begin{figure}
\centering
\includegraphics[width=0.75\textwidth]{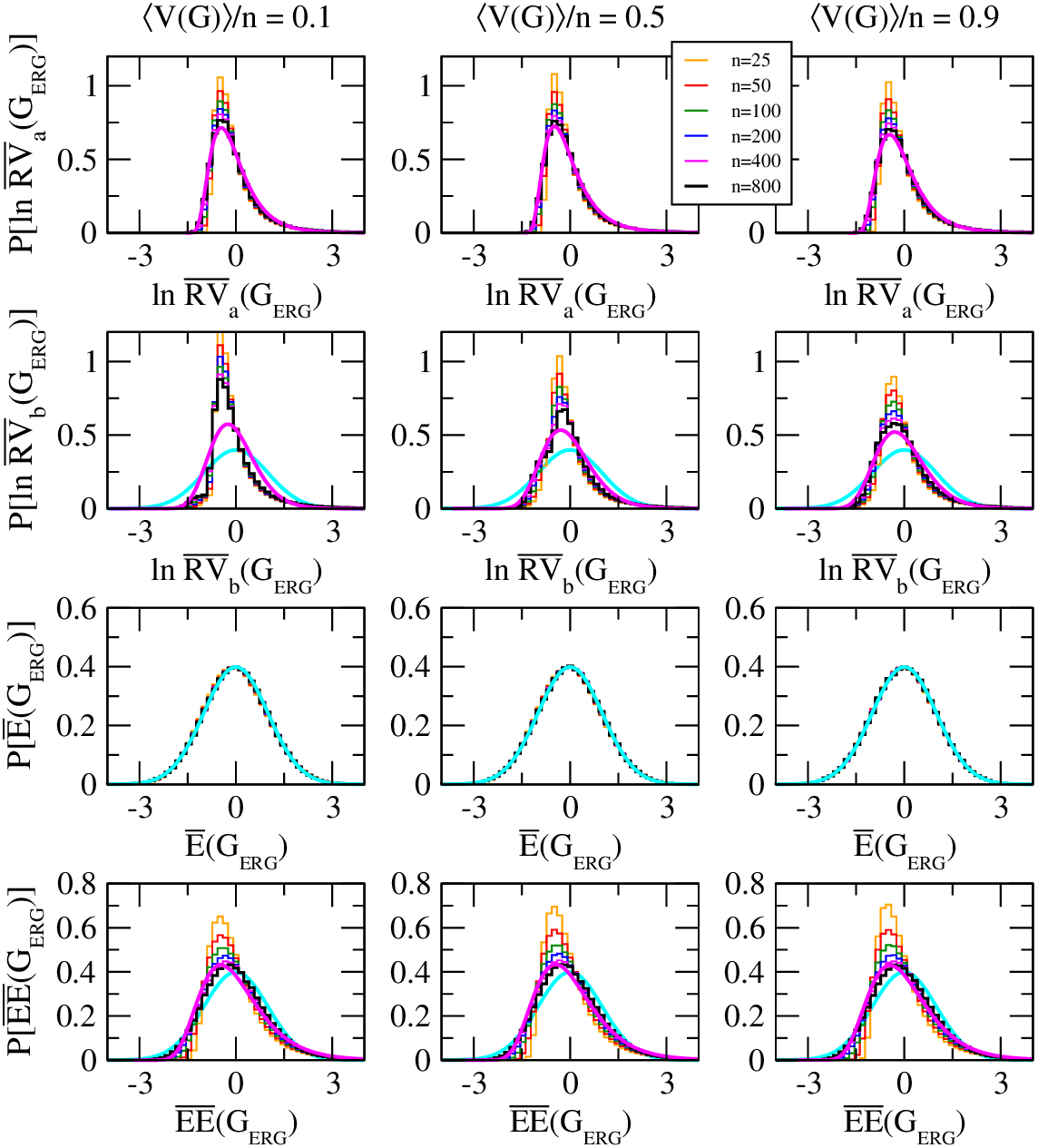}
\caption{Probability distribution functions of standardized spectral indices on Erd\H{o}s-R\'{e}nyi graphs: 
Rodr\'iguez-Vel\'azquez indices $\overline{RV_{a}}(G)$ and $\overline{RV_{b}}(G)$, graph energy 
$\overline{E}(G)$, and subgraph centrality $\overline{EE}(G)$.
Each panel displays six histograms corresponding to graphs of different sizes $n \in [25, 800]$.
Each column corresponds to a fixed value of the ratio $\bra V(G) \ket/n$.
Each histogram is constructed with $10^{6}$ values. 
Magenta full lines are fittings of Eq.~(\ref{ajuslogn}) to the distributions
corresponding to $n=800$; the values of the fitting parameters are reported in Table~\ref{ajustesER}.
The cyan full line in all panels is a normal distribution with zero mean and unit variance.}
\label{Fig11}
\end{figure}

In Fig.~\ref{Fig11} we present the probability distribution functions of standardized spectral indices on 
ERGs. As well as for MTIs, since RV indices and the exponential subgraph centrality grow exponentially 
with $\bra k \ket$, see e.g.~\cite{AMRS21}, we work with the distribution of their logarithms.

From Fig.~\ref{Fig11} we can clearly see that, remarkably, the shape of the distributions of spectral indices 
do not change with the graph connectivity (except for the distribution of $\ln RV_b(G)$ which shows a slight
dependence with $\bra V(G) \ket/n$); this in contrast with the distributions of degree-based indices whose
shapes evolve with $\bra V(G) \ket/n$. We also observe a slight dependence of the distribution
shapes with the graph size.

Moreover, note that only the distribution of $E(G)$ exhibits the shape of a normal distribution.
In contrast, the distributions of the Rodr\'iguez-Vel\'azquez indices as well as those of $EE(G)$ display an 
asymmetric, right-skewed, shape.
We found that the log-normal distribution function
\begin{equation}
    f(x, \sigma, \mu) = \frac{1}{ x  \sigma  \sqrt{2 \pi}}  \exp{\left(-\frac{(\log(x) - \mu)^2}{ 2 \sigma^2}\right)}-\beta,
    \label{ajuslogn}
\end{equation}
fits relatively well the distributions of the Rodr\'iguez-Vel\'azquez indices as well as those of $EE(G)$; see
the magenta lines in Fig.~\ref{Fig11} which are the fittings of Eq.~(\ref{ajuslogn}) to the distributions
corresponding to $n=800$.
Here, $\mu$ and $\sigma$ are the mean and standard deviation of the spectral indices, in a logarithmic 
scale, and $\beta$ is the distribution displacement on the $x$-axis.
The values of the fitting parameters are reported in Table~\ref{ajustesER}.

\begin{table}
\begin{center}
\begin{tabular}{| c | c | c | c |c | }
\hline
Index & $\bra V(G) \ket/n $ & $\sigma$ & $\mu$ & $\beta$\\ \hline
      & 0.1 & 0.4477 & 0.3206 &1.5821 \\
$RV_a$& 0.5 & 0.4918 & 0.2336&1.4781 \\
      & 0.9 & 0.4307 & 0.4183 &1.7211\\ \hline
      & 0.1 & 0.1857 & 1.3362 &3.9269 \\
$RV_b$& 0.5 & 0.2174 & 1.2587 &3.6506 \\
      & 0.9 & 0.2176 & 1.2814&3.7300 \\ \hline
      & 0.1 & 0.3252 & 1.1020 &3.1740\\
$EE$& 0.5 & 0.3146 & 1.1243 &3.2365 \\
    & 0.9 & 0.3095 & 1.1486 &3.3103 \\ \hline
\end{tabular}
\caption{Values of the parameters $\sigma$, $\mu$, and $\beta$ obtained from the
fittings of Eq.~(\ref{ajuslogn}) to the probability distribution functions (with $n=800$) of the spectral 
indices in Fig.~\ref{Fig11}.}
\label{ajustesER}
\end{center}
\end{table}

\section{Summary}
\label{conclusions}

In this work, we performed a thorough statistical (numerical) study of the probability distribution 
functions of topological and spectral indices on random graphs.
Specifically, we computed degree-based topological indices (TIs), degree-based multiplicative indices (MTIs)
Revan-degree indices, and spectral indices on two types of random graphs: Erd\"os-R\'enyi graphs (ERGs) and 
random geometric graphs (RGGs).

We performed our study by the use of the following indices.
\begin{itemize}
\item TIs: Number of non-isolated vertices, first and second Zagreb indices, Sombor index, Randi\'c index, and 
harmonic index.
\item MTIs: Narumi-Katayama index, multiplicative Zagreb indices, multiplicative Randi\'c index, and multiplicative 
harmonic index.
\item Revan-degree indices: First Revan Zagreb index, second Revan Zagreb index, Revan Sombor index, Revan 
Randi\'c index, and Revan harmonic index. Also, Revan Narumi-Katayama index, multiplicative Revan Zagreb 
indices, multiplicative Revan Randi\'c index, and multiplicative Revan harmonic index.
\item Spectral indice: Rodr\'iguez-Vel\'azquez indices, graph energy, and subgraph centrality.
\end{itemize}

It is relevant to mention that previos studies of the distributions of TIs were reported in 
Refs.~\cite{MMRS20,MMRS21,Y24,Y25}. However, the statistical studies of Refs.~\cite{MMRS20,MMRS21}
were not exhaustive while the analytical studies of Refs.~\cite{Y24,Y25} were not numerically verified, so
in this work we believe we fill those gaps.
Therefore, our results can be summarized as follows. 

For both ERGs and RGGs:
\begin{itemize}
\item[{\bf (i)}] asymptotically, for large enough connectivity and graph size, the distributions of the TIs follow a 
normal distribution (see Figs.~\ref{Fig03} and~\ref{Fig05});
 
\item[{\bf (ii)}]  asymptotically, for large enough connectivity and graph size, since the distributions of the logarithm 
of the MTIs follow a normal distribution (see Figs.~\ref{Fig04} and~\ref{Fig06}), the distributions of the MTIs follow 
a log-normal distribution;

\item[{\bf (iii)}] deep in the dense limit, the distributions of the Revan-degree indices follow a normal 
distribution (see Figs.~\ref{Fig07} and~\ref{Fig09});

\item[{\bf (iv)}] deep in the dense limit, since the distributions of the logarithm of the multiplicative 
Revan-degree indices follow a normal distribution (see Figs.~\ref{Fig08} and~\ref{Fig10}), the distributions of the 
multiplicative Revan-degree indices follow a log-normal distribution;

\item[{\bf (v)}] the distribution of the graph energy exhibits the shape of a normal distribution for any
graph connectivity and graph size (see Figs.~\ref{Fig11} and~\ref{Fig12});

\item[{\bf (vi)}] the distributions of the Rodr\'iguez-Vel\'azquez indices as well as those of the subgraph
centrality follow an asymmetric, right-skewed, log-normal distribution (see Eq.~(\ref{ajuslogn}) and Figs.~\ref{Fig11} 
and~\ref{Fig12}).
\end{itemize}

We finally stress that our results validate the analytical prediction of Ref.~\cite{Y23} stating that the distributions 
of TIs on ERGs converge to normal distributions; see Fig.~\ref{Fig03}.
However, and even more interesting, our results contradict the prediction of Ref.~\cite{Y25} which claims that
the limiting distribution of the Randi\'c index on RGGs is not the standard normal distribution; see Fig.~\ref{Fig05}.

We hope our results may motivate further analytical studies.

\acknowledgment{J.A.M.-B. thanks support from VIEP-BUAP (Grant No.~100405811-VIEP2025), Mexico.
C.T.M.-M. Thanks for the support from CONAHCYT (CVU No.~784756).}

% References should be ordered alphabetically.
\singlespacing

\newpage
\section*{Appendix: Distribution of topological and spectral indices on random geometric graphs}\label{appendix}
\addcontentsline{toc}{section}{Appendix: Random Geometric Graphs}
\label{appendix}

\begin{figure}
\centering
\includegraphics[width=0.75\textwidth]{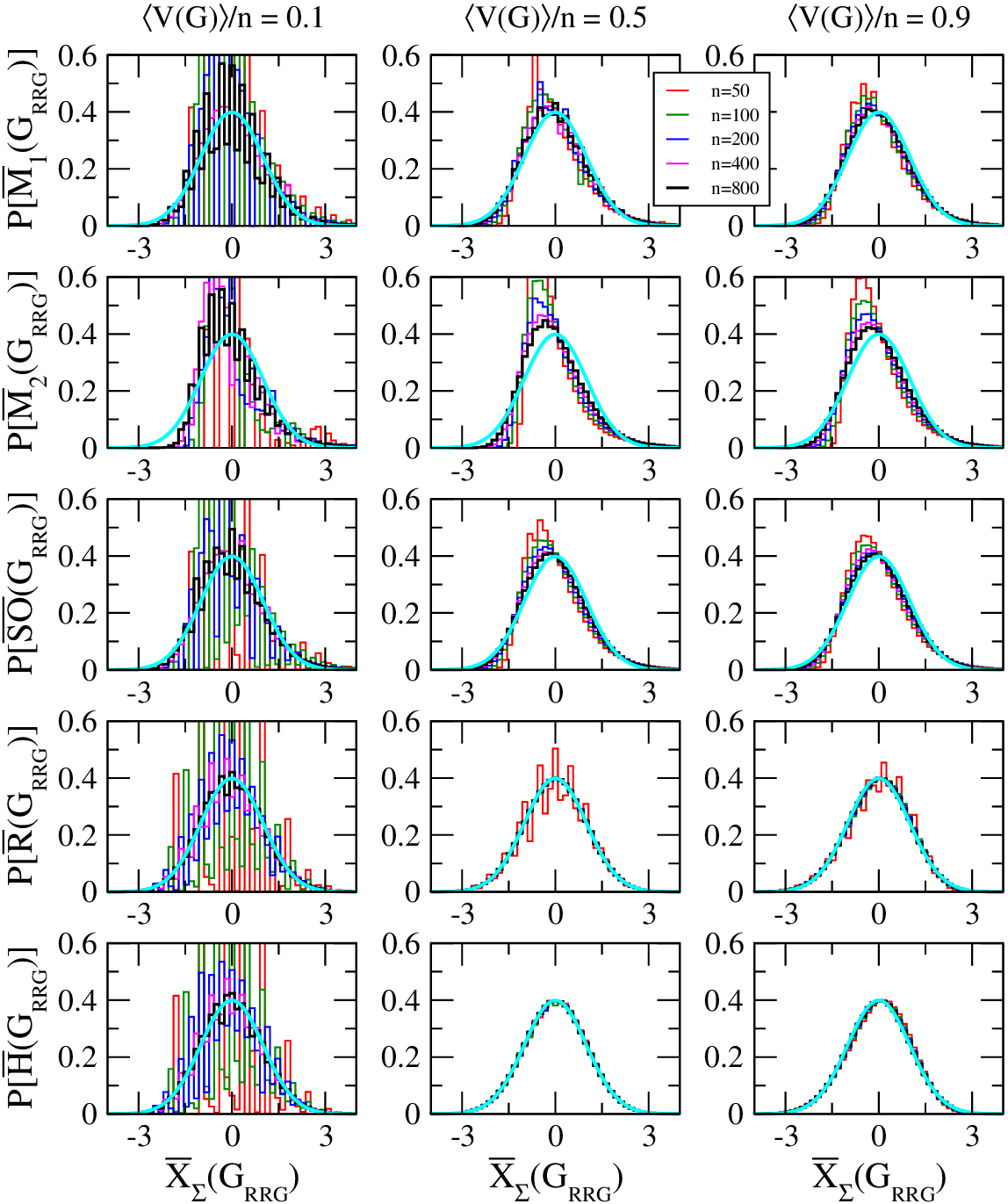}
\caption{Probability distribution functions of standardized degree-based topological indices
on random geometric graphs: First Zagreb index $\overline{M}_1(G)$, second Zagreb index 
$\overline{M}_2(G)$, Sombor index $\overline{SO}(G)$, Randi\'c index $\overline{R}(G)$, 
and harmonic index $\overline{H}(G)$.
Each panel displays five histograms corresponding to graphs of different sizes $n \in [50, 800]$.
Each column corresponds to a fixed value of the ratio $\bra V(G) \ket/n$.
Each histogram is constructed with $10^{6}$ values of $X_\Sigma(G)$.
The cyan full line in all panels is a normal distribution with zero mean and unit variance.}
\label{Fig05}
\end{figure}

\begin{figure}
\centering
\includegraphics[width=0.75\textwidth]{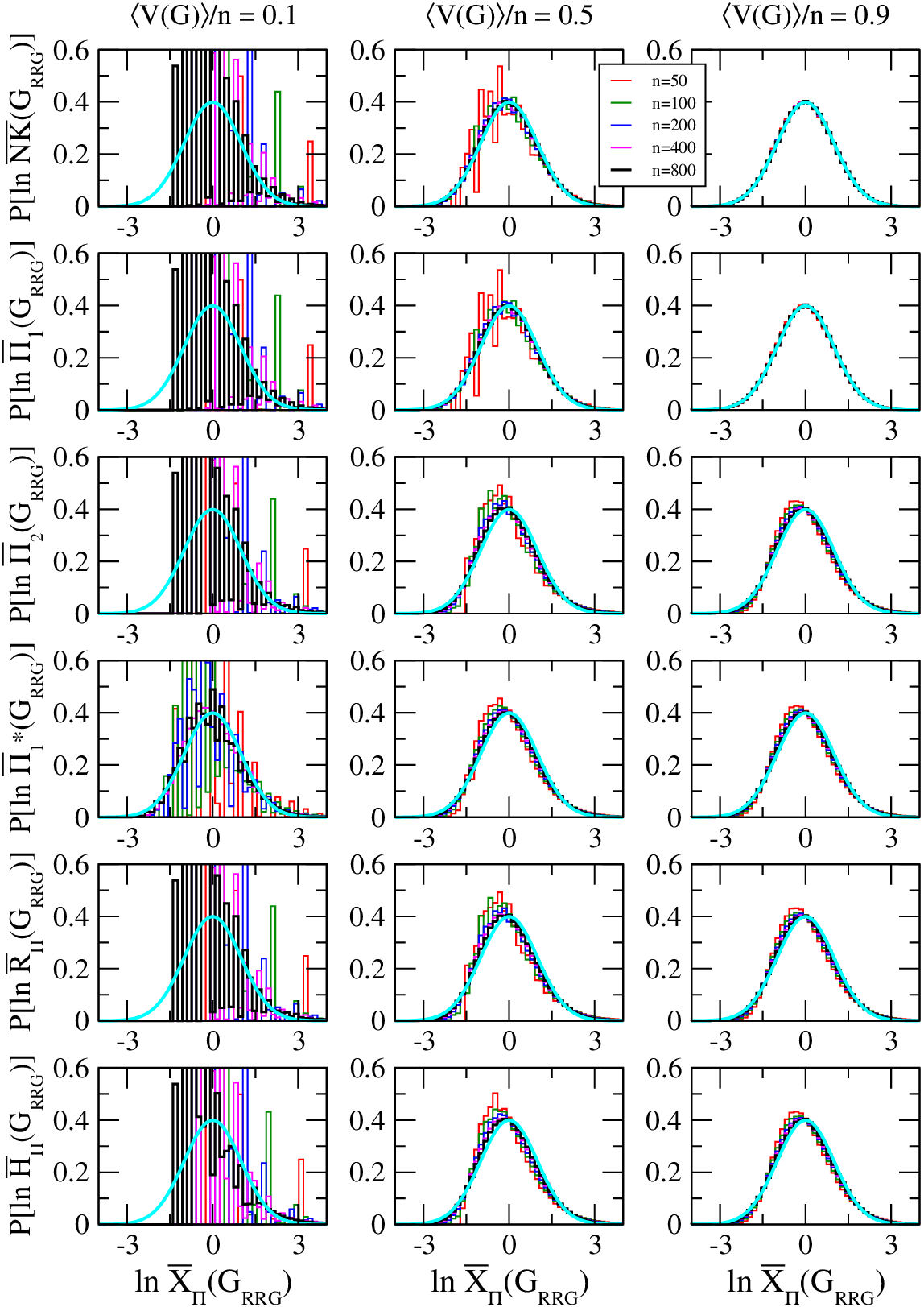}
\caption{Probability distribution functions of the logarithm of standardized 
multiplicative topological indices on random geometric graphs:
Narumi-Katayama index $\overline{NK}(G)$, multiplicative Zagreb 
indices $\overline{\Pi}_1(G)$, $\overline{\Pi}_2(G)$ and $\overline{\Pi}_1^*(G)$,
multiplicative Randi\'c index $\overline{R}_\Pi(G)$, and 
multiplicative harmonic index $\overline{H}_\Pi(G)$.
Each panel displays five histograms corresponding to graphs of different sizes $n$, ($n \in [50, 800]$).
Each column corresponds to a fixed value of the ratio $\bra V(G) \ket/n$.
Each histogram is constructed with $10^{6}$ values of $X_\Pi(G)$.
The cyan full line in all panels is a normal distribution with zero mean and unit variance.}
\label{Fig06}
\end{figure}

\begin{figure}
\centering
\includegraphics[width=0.75\textwidth]{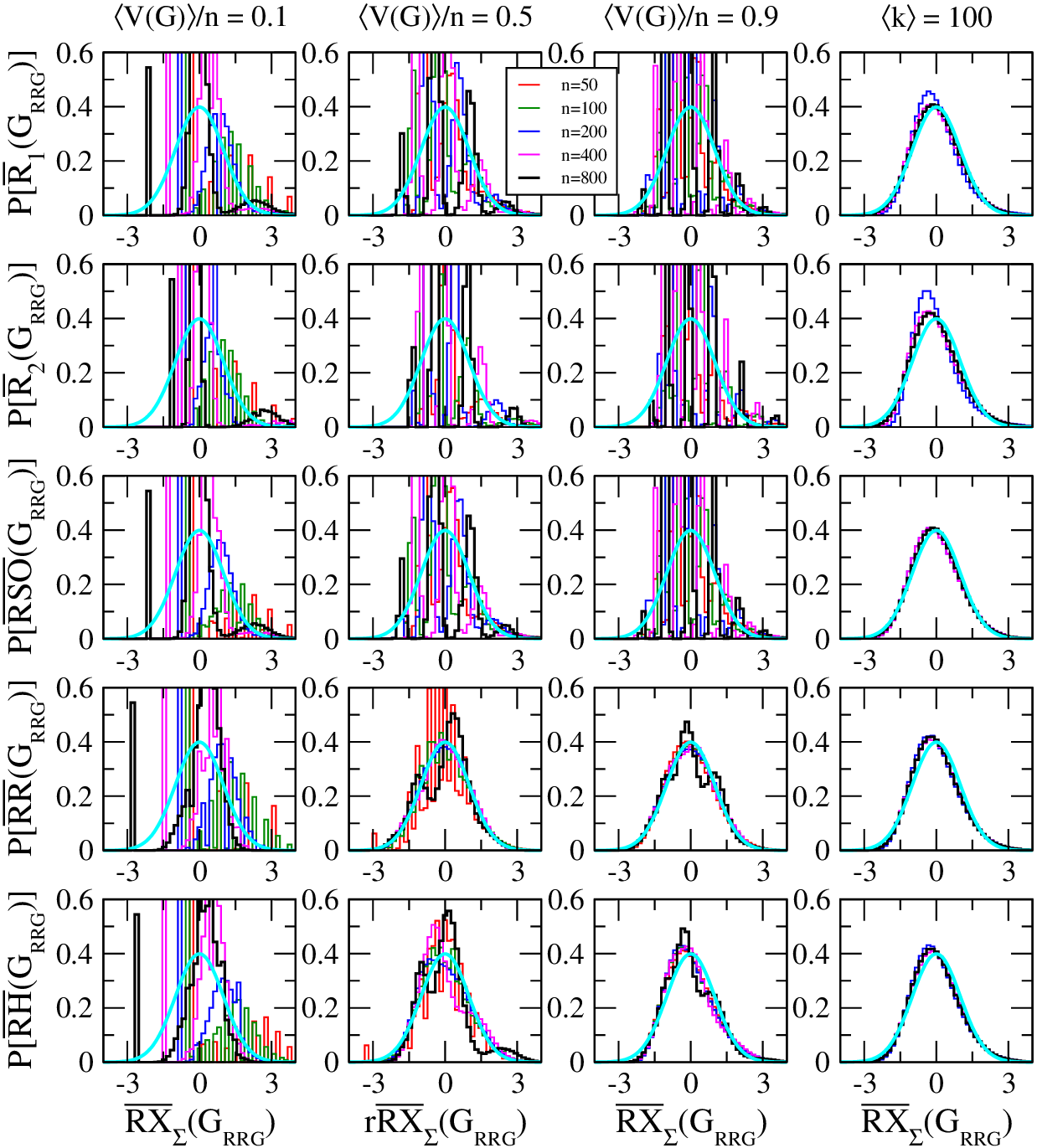}
\caption{Probability distribution functions of standardized Revan-degree indices on random geometric graphs: 
first Revan Zagreb index $\overline{R}_1(G)$, second Revan Zagreb index $\overline{R}_2(G)$, 
Revan Sombor index $\overline{RSO}(G)$, Revan Randi\'c index $\overline{RR}(G)$, and 
Revan harmonic index $\overline{RH}(G)$.
Each panel displays five histograms corresponding to graphs of different sizes $n \in [50, 800]$.
Each column corresponds to a fixed value of the ratio $\bra V(G) \ket/n$, except for the right 
column where $\bra k \ket=100$ is set.
Each histogram is constructed with $10^{6}$ values of $RX_\Sigma(G)$. 
The cyan full line in all panels is a normal distribution with zero mean and unit variance.}
\label{Fig09}
\end{figure}

\begin{figure}
\centering
\includegraphics[width=0.75\textwidth]{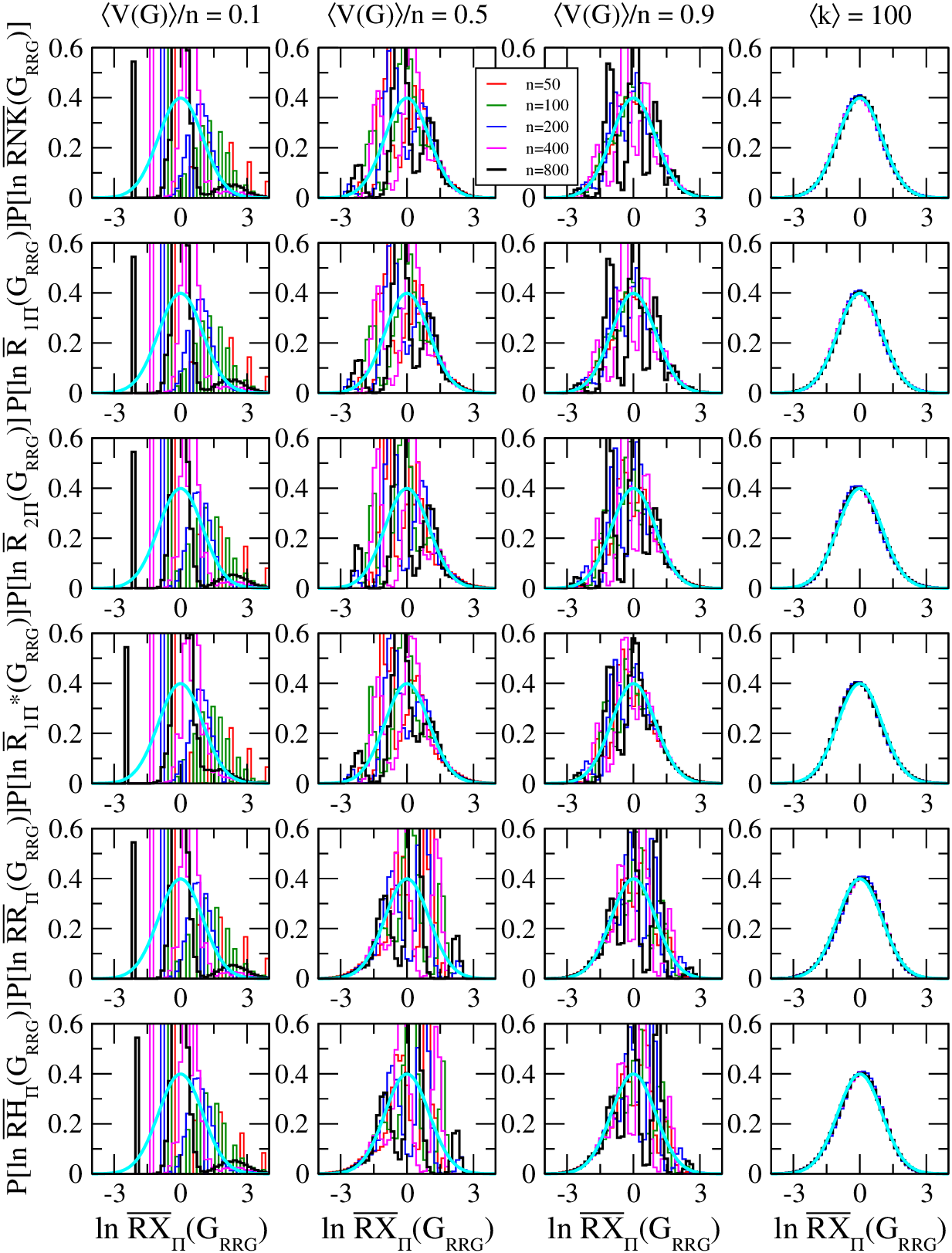}
\caption{Probability distribution functions of the logarithm of standardized multiplicative 
Revan-degree indices on random geometric graphs: 
Revan Narumi-Katayama index $\overline{RNK}(G)$, multiplicative Revan Zagreb indices 
$\overline{R}_{1\Pi}(G)$, $\overline{R}_{2\Pi}(G)$ and $\overline{R}_{1\Pi^*}(G)$,
multiplicative Revan Randi\'c index $\overline{RR}_\Pi(G)$, and multiplicative Revan 
harmonic index $\overline{RH}_\Pi(G)$.
Each panel displays five histograms corresponding to graphs of different sizes $n \in [50, 800]$.
Each column corresponds to a fixed value of the ratio $\bra V(G) \ket/n$, except for the right 
column where $\bra k \ket=10$ is set.
Each histogram is constructed with $10^{6}$ values of $RX_\Pi(G)$. 
The cyan full line in all panels is a normal distribution with zero mean and unit variance.}
\label{Fig10}
\end{figure}

\begin{figure}
\centering
\includegraphics[width=0.75\textwidth]{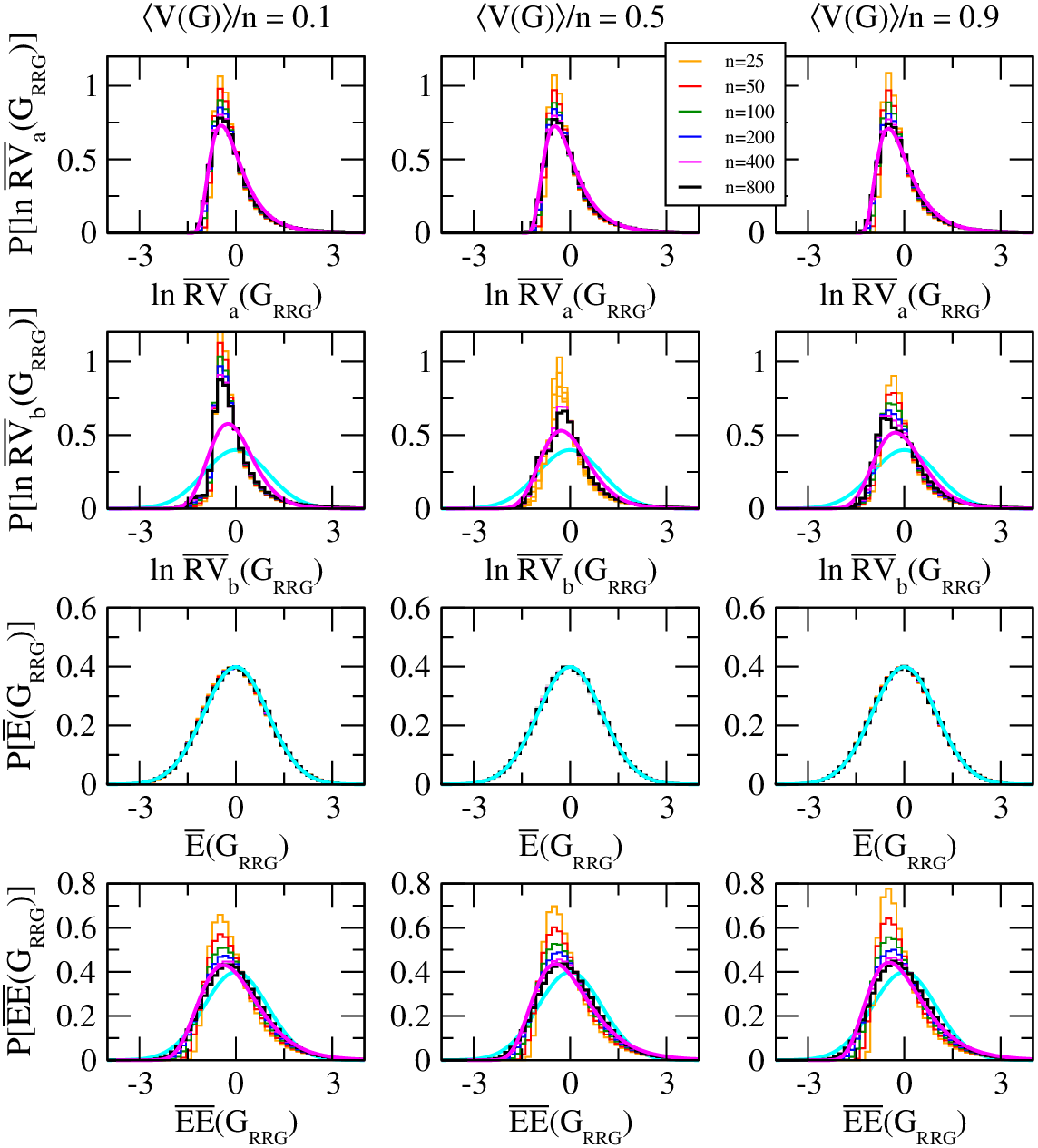}
\caption{Probability distribution functions of standardized spectral indices on
random geometric graphs: Rodr\'iguez-Vel\'azquez indices $\overline{RV_{a}}(G)$ and 
$\overline{RV_{b}}(G)$, graph energy $\overline{E}(G)$, and subgraph centrality $\overline{EE}(G)$.
Each panel displays six histograms corresponding to graphs of different sizes $n \in [25, 800]$.
Each column corresponds to a fixed value of the ratio $\bra V(G) \ket/n$.
Each histogram is constructed with $10^{6}$ values. 
Magenta full lines are fittings of Eq.~(\ref{ajuslogn}) to the distributions
corresponding to $n=800$; the values of the fitting parameters are reported in Table~\ref{ajustesRGG}.
The cyan full line in all panels is a normal distribution with zero mean and unit variance.}
\label{Fig12}
\end{figure}

In this Appendix, we report our results on the distributions of topological and spectral indices on RGGs, see
Figs.~\ref{Fig05}--\ref{Fig12}.
Note that Figs.~\ref{Fig05}--\ref{Fig12} for RGGs are equivalent to Figs.~\ref{Fig03}--\ref{Fig11} for ERGs, 
respectively. 
Indeed, from Figs.~\ref{Fig05}--\ref{Fig12} we draw similar conclusions as those already discussed in the 
main text for ERGs:
\begin{itemize}
\item[{\bf (i)}] Asymptotically, for large enough connectivity and graph size, the distributions of the TIs follow a normal 
distribution; see Fig.~\ref{Fig05}.

\item[{\bf (ii)}]  Asymptotically, for large enough connectivity and graph size, since the distributions of the logarithm of the 
MTIs follow a normal distribution (see Fig.~\ref{Fig06}), the distributions of the MTIs follow a log-normal distribution.

\item[{\bf (iii)}] Deep in the dense limit, the distributions of the Revan-degree indices follow a normal 
distribution; see Fig.~\ref{Fig09}.

\item[{\bf (iv)}] Deep in the dense limit, since the distributions of the logarithm of the multiplicative 
Revan-degree indices follow a normal distribution (see Fig.~\ref{Fig10}), the distributions of the multiplicative 
Revan-degree indices follow a log-normal distribution.

\item[{\bf (v)}] The distribution of the graph energy exhibits the shape of a normal distribution for any
graph connectivity and graph size; see Fig.~\ref{Fig12}.

\item[{\bf (vi)}] The distributions of the Rodr\'iguez-Vel\'azquez indices as well as those of the subgraph
centrality follow an asymmetric, right-skewed, log-normal distribution; see Fig.~\ref{Fig12} and Eq.~(\ref{ajuslogn}).

\end{itemize}

\begin{table}
\begin{center}
\begin{tabular}{| c | c | c | c |c | }
\hline
Index & $\bra V(G) \ket/n $ & $\sigma$ & $\mu$ & $\beta$\\ \hline
      & 0.1 & 0.4613 & 0.2753 &1.5226 \\
$RV_a$& 0.5 & 0.4763 & 0.2569&1.5029 \\
      & 0.9 & 0.4901 & 0.2581 &1.5114\\ \hline
      & 0.1 & 0.1838 & 1.3414 &3.9450\\
$RV_b$& 0.5 & 0.2067 & 1.3149 &3.8484 \\
      & 0.9 & 0.2032 & 1.3552&3.9998 \\ \hline
      & 0.1 & 0.2682 & 1.2775 &3.7256\\
$EE$& 0.5 & 0.3186 & 1.1073 &3.1873 \\
    & 0.9 & 0.3374 & 1.0406 &3.0007 \\ \hline
\end{tabular}
\caption{Values of the parameters $\sigma$, $\mu$, and $\beta$ obtained from the
fittings of Eq.~(\ref{ajuslogn}) to the probability distribution functions (with $n=800$) of the spectral 
indices in Fig.~\ref{Fig12}.}
\label{ajustesRGG}
\end{center}
\end{table}


\begin{thebibliography}{00}

\bibitem{R97}
M. Randi\'c, 
On characterization of chemical structure. 
{\it J. Chem. Inf. Comput. Sci.\/}, \textbf{37}(4) (1997) 672-687.

\bibitem{CT08}
V. Consonni, R. Todeschini,
New spectral indices for molecule description. 
{\it MATCH Commun. Math. Comput. Chem.} {\bf 60} (2008) 3--14.

\bibitem{T18}
N. Trinajstic, 
Chemical graph theory. CRC press.(2018).

\bibitem{G13}
I. Gutman,
Degree-based topological indices,
{\it Croat. Chem. Acta} {\bf 86} (2013) 351--361 .

\bibitem{R75}
Randi\'c, M.
On characterization of molecular branching.
{\it J. Am. Chem. Soc.} {\bf 97} (1975)  6609--6615.

\bibitem{G72}
I. Gutman and N. Trinajsti\'c,
Graph theory and molecular orbitals. Total $\phi$-electron energy of alternant hydrocarbons
{\it Chem. Phys. Lett.} {\bf 17} (1972) 535--538.

\bibitem{G21} 
I. Gutman, 
Geometric approach to degree--based topological indices: Sombor indices,
\textit{MATCH Commun. Math. Comput. Chem.\/} \textbf{86} (2021) 11--16.

\bibitem{NK84}
H. Narumi and M. Katayama,
Simple topological index – A newly devised index characterizing the topological nature of structural isomers of saturated hydrocarbons,
{\it Mem. Fac. Engin. Hokkaido Univ.} {\bf 16} (1984) 209-214.

\bibitem{TC10}
R. Todeschini and V. Consonni,
New Local Vertex Invariants and Molecular Descriptors
Based on Functions of the Vertex Degrees,
{\it MATCH Commun. Math. Comput. Chem.} {\bf 64} (2010) 359--372.

\bibitem{K17}
 V. R. Kulli, 
Revan indices of oxide and honeycomb graphs, 
\textit{Inter. J. Math. Appl.\/} \textbf{5}
(2017) 663--667.

\bibitem{K18}
 V. R. Kulli, 
F-Revan index and F-Revan polynomial of some families of benzenoid systems, 
\textit{J. Global Res. Math. Archives\/} \textbf{5} (2018) 1--6.

\bibitem{KG22}
 V. R. Kulli, I. Gutman, 
Revan Sombor index, 
\textit{J. Math. Inform.\/} \textbf{22} (2022) 23--27.

\bibitem{KMRS22}
V. R. Kulli, J. A. M\'endez-Bermudez, J. M. Rodriguez, J. M. Sigarreta, 
Revan Sombor indices: Analytical and statistical study,
\textit{Math. Biosciences Eng.\/} \textbf{20} (2023) 1801--1819.

\bibitem{AHMS22}
R. Aguilar-Sanchez, I. F. Herrera-Gonzalez, J. A. Mendez-Bermudez, and Jos\'e M. Sigarreta, 
Revan-degree indices on random graphs,
preprint arXiv:2210.04749.

\bibitem{RB19}
J.A. Rodriguez-Velazquez, A.T. Balaban, 
Two new topological indices based on graph adjacency matrix eigenvalues and eigenvectors,
\textit{J. Math. Chem.\/} \textbf{57} (2019) 1053--1074.

\bibitem{AMRS21}
R Aguilar-Sanchez, JA Mendez-Bermudez, Jos\'e M Rodr\'iguez, Jose M Sigarreta,
Analytical and statistical studies of Rodriguez-Velazquez indices,
{\it J. Math. Chem. } {\bf 59} (2021) 1246--1259.

\bibitem{LSG12}
X. Li, Y. Shi, and I. Gutman, 
Graph Energy, Springer, Berlin, 2012.

\bibitem{GR20}
I. Gutman and H. Ramane, 
Research on Graph Energies in 2019,
{\it MATCH Commun. Math. Comput. Chem.} {\bf 84}, (2020) 277--292.

\bibitem{ER05}
E. Estrada and J.A. Rodriguez-Velazquez, 
Subgraph centrality in complex networks,
{\it Phys. Rev. E} {\bf 71}, (2005) 056103.

\bibitem{MMRS20}
C. T. Mart\'{\i}nez-Mart\'{\i}nez, J. A. Mendez-Ber\'mudez, J. M. Rodr\'{\i}guez, and J. M. Sigarreta,
Computational and analytical studies of the Randi\'c index in Erd\"os--R\'enyi models,
{\it Appl. Math. Comput.} {\bf 377}, (2020) 125137.

\bibitem{MMRS21}
C. T. Mart\'{\i}nez-Mart\'{\i}nez, J. A. Mendez-Bermudez, J. M. Rodr\'{\i}guez, and J. M. Sigarreta,
Computational and analytical studies of the harmonic index in Erd\"os--R\'enyi models,
{\it MATCH Commun. Math. Comput. Chem.} {\bf 85},  (2021) 395--426.

\bibitem{AMRS20}
R. Aguilar-Sanchez, J. A. Mendez-Bermudez, F. A. Rodrigues, and J. M. Sigarreta-Almira,
Topological versus spectral properties of random geometric graphs,
{\it Phys. Rev. E} {\bf 102}, (2020) 042306.

\bibitem{AHMS20}
R. Aguilar-Sanchez, I. F. Herrera-Gonzalez, J. A. Mendez-Bermudez, and J. M. Sigarreta,
Computational properties of general indices on random graphs.
{\it Symmetry} {\bf 12}, (2020) 1341.

\bibitem{AMRS24}
R. Aguilar-Sanchez, J. A. M\'endez-Berm\'udez, J. M. Rodriguez, J. M. Sigarreta,
Multiplicative topological indices: Analytical properties and application to random graphs,
{\it AIMS Mathematics} {\bf 9}, (2024) 3646.

\bibitem{Y24}
M. Yuan,
Asymptotic distribution of degree-based topological indices,
{\it MATCH Commun. Math. Comput. Chem.} {\bf 91}, (2025) 135--196.

\bibitem{Y25}
M. Yuan,
Limiting distribution for the Randic index of a random geometric graph,
{\it MATCH Commun. Math. Comput. Chem.} {\bf 93}, (2024) 767--789.

\bibitem{Y23}
M. Yuan,
On the Randic index and its variants of network data, 
TEST {\bf 33} (2023) 155--179.

\bibitem{YZ23}
M. Yuan and X. Zhao,
Asymptotic distributions of the average clustering coefficient and its variant,
preprint arXiv:2311.10979.

\bibitem{SR51}
R. Solomonoff and A. Rapoport,
Connectivity of random nets.
{\it  Bull. Math. Biophys.} {\bf 13}, (1951) 107--117.

\bibitem{ER59}
P. Erd\"os and A. R\'enyi,
On random graphs.
{\it  Publ. Math. (Debrecen)} {\bf 6}, (1959) 290--297.

\bibitem{ER60}
P. Erd\"os and A. R\'enyi,
On the evolution of random graphs,
Inst. of the Hung. Acad. of Sci. {\bf 5}, (1960) 17--61;
P. Erd\"os and A. R\'enyi,
On the strength of connectedness of a random graph,
{\it  Acta Mathematica Hungarica} {\bf 12}, (1961) 261--267 .

\bibitem{DC02}
J. Dall and M. Christensen,
Random geometric graphs,
{\it  Phys. Rev. E} {\bf 66}, (2002) 016121.

\bibitem{P03}
M. Penrose,
Random Geometric Graphs, Oxford University Press, Oxford, 2003.

\bibitem{EM15}
E. Estrada and M. Sheerin,
Random rectangular graphs,
{\it  Phys Rev. E} {\bf 91}, (2015) 042805.

\bibitem{F87}
S. Fajtlowicz,
On conjectures of Graffiti--II,
{\it Congr. Numer.} {\bf 60} (1987) 187--197.



%\bibitem{AMRS20}
%R. Aguilar-Sanchez, J. A. Mendez-Bermudez, F. A. Rodrigues, and J. M. Sigarreta-Almira,
%Topological versus spectral properties of random geometric graphs,
%Phys. Rev. E 102, 042306 (2020).
%
%\bibitem{AMRS21}
%R. Aguilar-Sanchez, J. A. Mendez-Bermudez, J. M. Rodriguez, and J. M. Sigarreta,
%Normalized Sombor indices as complexity measures of random graphs,
%Entropy 23, 976 (2021).



%\bibitem{BQRS} P. Bosch, Y. Quintana, J. M. Rodr\'{\i}guez, J. M. Sigarreta,
%Jensen-type inequalities for m-convex functions, to appear in {\it Open Math.}
%
%\bibitem{GutmanMilovanovic} I. Gutman, I. Milovanovi\'c, E. Milovanovi\'c,
%Relations between Ordinary and Multiplicative Degree-Based Topological Indices,
%{\it Filomat} {\bf 32:8} (2018) 3031--3042.
%
%\bibitem{Kober} H. Kober, On the arithmetic and geometric means and on H\"older's inequality, {\it Proc. Amer. Math. Soc.} {\bf 9} (1958) 452--459.
%
%\bibitem{Petrovic} M. Petrovi\'c's, Sur Une Fonctionnelle, {\it Publ. Math. Univ. Belgrade} {\bf 1} (1932) 146--149.
%
%\bibitem{RetiGutman} T. R\'eti, I. Gutman, Relations between ordinary and multiplicative Zagreb indices,
%{\it Bull. Inter. Math. Virtual Inst.} {\bf 2} (2012) 133--140.
%
%\bibitem{ZGA} B. Zhou, I. Gutman, T. Aleksi\'c, A note on Laplacian energy of graphs, {\it MATCH Commun. Math. Comput. Chem.} {\bf 60} (2008) 441--446.








\end{thebibliography}
\end{document}